\documentclass[paper=a4,english,fontsize=11pt,parskip=half,abstract=true]{scrartcl}
\usepackage{babel}
\usepackage[utf8]{inputenc}
\usepackage[T1]{fontenc}
\usepackage[left=20mm,right=20mm,top=30mm,bottom=30mm]{geometry}
\usepackage{amsmath}
\usepackage{amsthm}
\usepackage{amssymb}
\usepackage{enumerate} 
\usepackage{mathtools}
\mathtoolsset{centercolon} 
\usepackage[bookmarks=true,
            pdftitle={A counterexample to a conjecture of A. R. Miller},
            pdfauthor={Gabriel Navarro, Benjamin Sambale},
            pdfkeywords={},
            pdfstartview={FitH}]{hyperref}

\newcommand{\RR}{\mathbb{R}}
\newcommand{\Irr}{\operatorname{Irr}}
\newcommand{\tr}{\operatorname{tr}}

\title{A counterexample to a conjecture\\ of A.~R.~Miller}
\author{Gabriel Navarro\footnote{Department of Mathematics, Universitat de Val\`encia, 46100 Burjassot,
        Spain, \href{mailto:gabriel@uv.es}{gabriel@uv.es}} \ and Benjamin Sambale\footnote{Institut für Algebra, Zahlentheorie und Diskrete Mathematik, Leibniz Universität Hannover, Welfengarten 1, 30167 Hannover, Germany,
\href{mailto:sambale@math.uni-hannover.de}{sambale@math.uni-hannover.de}}}
\date{\today}

\begin{document}
\frenchspacing
\maketitle
\begin{abstract}\noindent
Let $\chi$ be an irreducible character of a finite group $G$. 
A.~R.~Miller conjectured that the proportion of elements $g\in G$ such that $\chi(g)$ is zero or a root of unity is at least $1/2$. We construct a character of a perfect group of order $69{,}120$ such that this proportion is $511/1152\approx 0.44$. 
\end{abstract}


Let $x\in\RR\setminus\{1\}$ be an algebraic integer such that all Galois conjugates $x_1,\ldots,x_n$ of $x$ are real and positive. A well-known theorem of Siegel~\cite[Theorem~III]{Siegel} states that 
\[\overline{\tr}(x):=\frac{x_1+\ldots+x_n}{n}\ge\frac{3}{2}\]
with equality if and only if $x=\frac{3\pm\sqrt{5}}{2}$. 
Thompson has observed that this theorem applies to $x=|\chi(g)|^2$ for any character $\chi$ of finite groups $G$ and $g\in G$.
When $\chi$ is irreducible, the orthogonality relations imply that the proportion of elements $g\in G$ such that $\chi(g)$ is zero or a root of unity is at least $1/3$, i.\,e.
\[\lambda(\chi):=\bigl|\bigl\{g\in G:\chi(g)=0\text{ or }|\chi(g)|=1\bigr\}\bigr|\ge\frac{|G|}{3}\]
(see \cite[Problem~3.15]{Isaacs}).

Siegel's theorem has been improved many times, but the exact value of
\[\lambda:=\liminf_{x}\overline{\tr}(x)\]
is still unknown. Most recently, Orloski--Sardari--Smith~\cite[p. 2010]{SSSproblem} proved that $\lambda\ge 1.80203$. On the other hand, the sequence $\{4\cos(2\pi/p)^2:p\text{ prime}\}$ shows that $\lambda\le 2$. 

Since character values are cyclotomic integers (i.\,e. integers of a cyclotomic field), it is expected that stronger bounds can be given. In fact, Cassels~\cite[Lemma~3]{Cassels} proved for every cyclotomic integer $x$ that $\overline{\tr}(|x|^2)\ge 2$ unless $x$ is a sum of at most two roots of unity. 
Perhaps motivated by this result, Miller~\cite[Conjecture~1]{ARMiller} conjectured that $\lambda(\chi)\ge|G|/2$ for all $\chi\in\Irr(G)$ (see also \cite{ARMiller2,ARMiller3}). He proved this bound for all nilpotent groups, among other cases. This was extended to all groups with a Sylow tower by Moretó--Navarro~\cite[Theorem~C]{MoretoNavarroMiller}. The conjecture seems to be open for solvable groups. 

Using GAP~\cite{GAPnew}, we found a perfect group $G\cong C_{12}.(C_2^4\rtimes A_6)$ of order $69{,}120$ with eight (Galois conjugate) irreducible characters $\chi$ of degree $12$ such that $\lambda(\chi)=511/1152\approx 0.44$. 
This can be checked with the following code:
\begin{verbatim}
n:=69120;
G:=PerfectGroup(n,2);
ct:=CharacterTable(G);
chi:=First(Irr(ct), c -> c[1]=12);
pos:=PositionsProperty(chi, x -> x=0 or x*ComplexConjugate(x)=1);
Sum(SizesConjugacyClasses(ct){pos})/n;
\end{verbatim}

In accordance with Siegel's and Cassels's theorems mentioned above, there exists $g\in G$ such that $\chi(g)=\zeta_5+\overline{\zeta_5}=2\cos(2\pi/5)$ and $|\chi(g)|^2=\frac{3-\sqrt{5}}{2}$ in this example. In fact, there is no $g\in G$ such that $\chi(g)$ is a root of unity. There are more such examples among perfect groups involving $A_6$, but we did not find a smaller value of $\lambda(\chi)$.  
It remains an open problem to determine $\min\lambda(\chi)$ over all finite groups.

\section*{Acknowledgment}
The research of the first author is supported by grant PID2022-137612NB-I00 funded by MCIN/AEI/ 10.13039/501100011033 and ERDF “A way of making Europe”.
The second author thanks Marco Soriano for maintaining a remote machine at the LUH, on which we found this counterexample during the night.

\end{document}